\documentclass[12pt,a4paper]{article}

\usepackage{amsmath,amssymb,graphicx,epsfig}

\begin{document}

\title{`Strange' occurrences in {\em SuperEnalotto}}

\author{Germano D'Abramo\vspace{0.2cm}\\
{\small Istituto Nazionale di Astrofisica,}\\ 
{\small Via Fosso del Cavaliere 100,}\\
{\small 00133, Roma, Italy.}\\
{\small E--mail: {\tt Germano.Dabramo@iasf-roma.inaf.it}}}
\vspace{0.2cm}
\date{}

\maketitle

\begin{abstract}

In this paper a way is suggested for calculating the probability of {\em 
consecutive number} strings within a sequence of $n$ numbers randomly 
drawn (without replacement) among the set of the first $N$ consecutive 
numbers, with $N\gg n$.

An explicit derivation is carried out for the special case of {\em 
SuperEnalotto}, nowadays the most famous lottery in Italy, with $N=90$ 
and $n=6$. It turns out that, on average, one every three drawings 
presents one or more {\em consecutive number} strings inside.\\

\end{abstract}

\section{Introduction}

Among lotteries and gambling games endorsed by the Italian government, 
{\em SuperEnalotto} is surely the most popular one, even in the neighbor 
states like France, Austria and Switzerland. It has been introduced for 
the first time on December 3, 1997 and since then it has provided some 
of the biggest payoffs of all lotteries in the world~\cite{sel}.

The game is currently issued three times a week. The player gains the 
jackpot (which increases at every issue, depending on the total number of 
players each time and on the number of the previous not winning issues) 
if she matches the 6 numbers drawn out of 90 (from 1 to 90), regardless 
of the order in which the numbers are drawn.

Concerning the winning odds, {\em SuperEnalotto} is considered to be one 
of the most difficult gambling game among the existing ones: winning the 
jackpot means to match the drawn un-ordered collection of six numbers 
among all the possible 6-combinations of 90 numbers. It is easy to see 
that there are ${90 \choose 6}$ of these combinations, so the winning 
odds at every issue of the game are $1$ in $622,614,630$.

To give an idea of the tininess of the figures at play, the odds {\bf 
you} have to gain the jackpot in the next issue of the game are 
comparable~\cite{stu} with the odds {\bf you} have to be hit, here on 
the Earth, by a kilometre-sized asteroid during the very same day you 
play the game! Obviously, this does not mean that no one will never gain 
the jackpot: many people play the game, thus the expectation value of 
winnings is not tiny at all.

Although the author is not a greedy player of {\em 
SuperEnalotto}\footnote{Probably because he is pretty confident that he 
will never be hit by a kilometre-sized asteroid in his whole life.}, his 
attention was drawn by an apparently counter-intuitive high frequency 
with which some peculiar number patterns appear within the six drawn 
numbers. In particular, he noted an unexpected high 
occurrence\footnote{Indeed, this is a personal impression, not an 
objective fact.} of {\em consecutive numbers}, namely like those in the 
last drawings of year 2009 as reported in Table~\ref{tab1}.

Among conspiracy-inclined people, and there is plenty of them between 
greedy players, this feature is read as the unequivocal sign that 
drawings are intentionally biased by the lottery operator.

On the contrary, the first though of the author was to understand if the 
laws of probability are able to explain such an apparently bizarre 
behavior. Obviously, they do and in the following Section it is shown 
how.

\begin{table}[t]
\begin{center}
{\tiny 
\begin{tabular}{lcccccc}
\hline
\multicolumn{6}{c}{Date \qquad\qquad\qquad Drawn numbers}\\
\hline
        31/12/2009  &  {\bf 40} &  {\bf 41}  & 	45  & 	{\bf 51} &  {\bf 52} & 79 \\ 	       
      	29/12/2009  & 11  &	16  &	21  &	62  &	70  &	84 \\      
      	28/12/2009 &  34  &	36   &	{\bf 61}  & {\bf 62}   & 78   &	88  \\         
      	24/12/2009 &  5   & 25   &	37   &	52   &	62   &	79  \\        
      	22/12/2009 &  6   &	13  & 	17   &	52   &	64   &	82  \\         
      	19/12/2009 &  9   &	15   &	17   &	34   &	73   &	83  \\        
      	17/12/2009 &  1   &	23   &	32   &	53   &	66   &	83  \\        
      	15/12/2009 &  1   &	4   &	48   &	53   &	61   &	85  \\        
      	12/12/2009 &  10   &	16   &	40   &	55   &	66   &	77  \\        
      	10/12/2009 &  33   &	36   &	39   &	50   &	65   &	70  \\      
      	09/12/2009 &  9   &	30   &	50   &	57   &	62   &	70  \\        
      	05/12/2009 &  20   &	63   &	66   &	{\bf 74}   &	{\bf 75}   & 77  \\        
      	03/12/2009 &  19   &	30   &	34   &	49   &	54   &	74  \\       
      	01/12/2009 &  21   &	25   &	36   &	51   &	82   &	90  \\        
      	28/11/2009 &  3   &	42   &	52   &	67   &	78   &	86  \\         
      	26/11/2009 &  26   &	{\bf 31}   &	{\bf 32}   &	50   &	61   &	77  \\          
      	24/11/2009 &  19   &	47   &	69   &	71   &	{\bf 86}   &	{\bf 87}  \\         
      	21/11/2009 &  {\bf 36}   & {\bf 37}   &	50   &	57   &	73   &	83  \\        
      	19/11/2009 &  1   &	35   &	42   &	55   &	57   &	62  \\         
      	17/11/2009 &  16   &	25   &	28   &	30   &	39   &	51  \\        
      	14/11/2009 &  15   &	27   &	31   &	40   &	63   &	82  \\         
      	12/11/2009 &  37   &	43   &	49   &	66   &	73   &	88  \\          
      	10/11/2009 &  24   &	32   &	56   &	59   &	76   &	81  \\        
      	07/11/2009 &  9   &	11   &	70   &	{\bf 74}   & {\bf 75}   & 84  \\         
      	05/11/2009 &  12   &	20   &	28   &	57   &	74   &	85  \\         
      	03/11/2009 &  22   &	35   &	44   &	66   &	71   &	88  \\         
      	31/10/2009 &  13   &	16   &	41   &	49   &	58   &	76  \\          
      	29/10/2009 &  29   &	37   &	42   &	52   &	54   &	59  \\          
      	27/10/2009 &  38   &	{\bf 43}   &	{\bf 44}  &	{\bf 45}   &	63   &	77  \\          
      	24/10/2009 &  6   &	26   &	{\bf 37}   &	{\bf 38}   &	55   &	73  \\        
      	22/10/2009 &  29   &	{\bf 33}   &	{\bf 34}    &	{\bf 64}   & {\bf 65}   & 87  \\         
      	20/10/2009 &  41   &	47   &	49   &	55   &	77   &	90  \\           
      	17/10/2009 &  8   &	24   &	45   &	47   &	72   &	84  \\        
      	15/10/2009 &  3   &	29   &	51   &	55   &	73   &	81  \\        
      	13/10/2009 &  5   &	11   &	40   &	43   &	56   &	82  \\
      	10/10/2009  &  9   & 	{\bf 33}   & {\bf 34}   &	52   &	{\bf 82}   &	{\bf 83}  \\          
      	08/10/2009 &  1   &	12   &	17   &	65   &	70   &	90  \\          
      	06/10/2009 &  {\bf 23}   & {\bf 24}   &	49   &	{\bf 69}   &	{\bf 70}   &	72  \\          
      	03/10/2009 &  28   &	38   &	{\bf 59}   & {\bf 60}   &	71   &	83  \\          
      	01/10/2009 &  4   &	19   &	34   &	52   &	68   &	77   \\       
\hline
\end{tabular}
}
\end{center}
\label{tab1}
\caption{{\em SupeEnalotto} drawings in the last three months of year 2009.
In bold: {\em consecutive number} strings.}
\end{table}

\section{Consecutive number probability}

First of all, let us define the {\em consecutive number} probability 
${\cal C}(n,N)$: this is the probability to have one or more strings of 
consecutive numbers (of all lengths, starting from strings of 2 
consecutive numbers up to strings of $n$ consecutive numbers) within a 
sequence of $n$ numbers randomly drawn from the set $\{1,2,3,4,...,N\}$ 
of the first $N$ consecutive numbers. For the sake of simplicity, the 
{\em consecutive number} probability is explicitly derived here only 
for the special case with $N=90$, and $n=6$, namely the case of {\em 
SuperEnalotto}. Generalizations of such a result should not pose 
particular difficulties.

In what follows every drawn sequence is thought as written with numbers 
in increasing order (as in Table~\ref{tab1}) and in the following 
format:

\begin{equation}
s_1\,a_1\,s_2\,a_2\,s_3\,a_3\,s_4\,a_4\,s_5\,a_5\,s_6\,a_6\,s_7,
\label{eq1}
\end{equation}
where the $a_i$ are the drawn numbers (with $a_i<a_j$, if $i<j$), while 
the $s_k$ are numbers which express the distance between two contiguous 
$a_i$. Namely, $s_2=(a_2-a_1)-1$, $s_3=(a_3-a_2)-1$, and so on, while 
$s_1=a_1-1$ and $s_7=(90-a_6)$. As a particular case $s_1=0$, if 
$a_1=1$, and $s_7=0$, if $a_6=90$.

An important arithmetical relation then holds:

\begin{equation}
s_1+s_2+s_3+s_4+s_5+s_6+s_7=90-6=84,
\label{eq2}
\end{equation}
with $0\leq s_k\leq 84$.

Equation~(\ref{eq2}) allows to cope with consecutive numbers within the 
drawn sequence: as a matter of fact, one has consecutive number strings 
when at least one of the value $\{s_2,s_3,s_4,s_5,s_6\}$ is equal to 
zero. The values of $s_1$ and $s_7$ are not important for the 
consecutive number issue.

Thus, in order to calculate the numerical value of ${\cal C}(6,90)$, one 
has to count all the ways of writing the integer 84 as a sum of the 
numbers $\{s_1$, $s_2$, $s_3$, $s_4$, $s_5$, $s_6$, $s_7\}$, with at 
least one of the value $\{s_2,s_3,s_4,s_5,s_6\}$ equal to zero: in fact, 
${\cal C}(6,90)$ is equal to the ratio between this number of ways (let 
call it $N_c$) and the ways of writing 84 as the sum of the numbers 
$\{s_1$, $s_2$, $s_3$, $s_4$, $s_5$, $s_6$, $s_7\}$, without any proviso 
on their values, namely, with $0\leq s_k\leq 84$. Let call this last 
number $N_t$.

This is a partition problem and like most of partition problems it can 
be solved thanks to the technique of {\em generating functions}.

Consider first the calculation of the denominator of the above 
probability fraction, the number $N_t$.

Every $s_k$ spans from $0$ to $84$; thus, let us expand the following 
polynomial power:

\begin{equation}
(1+x+x^2+x^3+\cdots +x^{82}+x^{83}+x^{84})^7=\bigg(\frac{x^{85}-1}{x-1}\biggl)^7.
\label{eq3}
\end{equation}

After the expansion, the multiplicative coefficient of the term $x^{84}$ 
provides exactly the number $N_t$: such a coefficient sums up the number 
of way in which the term $x^{84}$ can be obtained as a product of the 
seven terms $x^l$ (one for each $s_k$), and thus, it sums up also the 
ways in which the exponent of $x^{84}$ can be obtained as a sum of the 
seven exponents $l$, with $0\leq l\leq 84$.

Applying multiple derivative to eq.~(\ref{eq3}) and making a suitable 
normalization (to simplify the numerical coefficient equal to $84!$ 
which originates from the exponent of $x^{84}$ after multiple derivative), 
a relatively simple algorithm for calculating the coefficient of $x^{84}$ 
can be obtained as follows:

\begin{equation}
N_t=\frac{1}{84!}\cdot\frac{d^{84}}{dx^{84}}
\biggl[\bigg(\frac{x^{85}-1}{x-1}\biggl)^7\biggr]
\Biggl|_{x=0}=622,614,630.
\label{eq4}
\end{equation}

The numerical result in eq.~(\ref{eq4}) has been obtained in less than a 
minute with free math software {\tt Sage}~\cite{sage}, running on a not 
so brand-new laptop. The reader must have recognized such a figure. It is 
the total number of 6-combinations of the set of 90 numbers, namely ${90 
\choose 6}$, as it should be.

The next step is to calculate the number $N_c$. This task is a bit more 
involved, although it does not require anything new with respect to the 
derivation of the number $N_t$. According to the format of the drawn 
sequence introduced in eq.~(\ref{eq1}), one has {\em consecutive 
numbers} when at least one of the five numbers $\{s_2$, $s_3$, $s_4$, 
$s_5$, $s_6\}$ is equal to zero.

Let us start with the simplest case when only one $s_k$ among $\{s_2$, 
$s_3$, $s_4$, $s_5$, $s_6\}$ is equal to zero while the other ones are 
different from zero. This case corresponds to the situation in which 
there is only a single pair of {\em consecutive numbers} within the 
drawn sequence.

First of all, there are ${5 \choose 1}=5$ cases in which one such $s_k$ 
may be equal to zero, since there are exactly five numbers $s_k$. 
Moreover, for each one of such cases one must count the ways in which 
the relation~(\ref{eq2}) is fulfilled with one $s_k$ equal to zero and 
the other ones assuming values $1\leq s_k\leq 84$.

The arguments introduced above suggest that the number of ways is 
equal to the coefficient of the term $x^{84}$ in the expansion of the 
following polynomial powers product:

\begin{equation}
(1+x+x^2+\cdots +x^{83}+x^{84})^{\bf 2}(x+x^2+\cdots +x^{83}+x^{84})^{\bf 4}
=\bigg(\frac{x^{85}-1}{x-1}\biggl)^{\bf 2}\bigg(\frac{x^{85}-1}{x-1}-1\biggl)^{\bf 4},
\label{eq5}
\end{equation}
where the factor to the left on both sides (that with exponent equal to 
${\bf 2}$) takes into account the behavior of terms $s_1$ and $s_7$ 
(that can be also equal to zero), while the factor to the right on 
both sides (that with exponent equal to ${\bf 4}$) takes into account 
the behavior of the four terms among $\{s_2$, $s_3$, $s_4$, $s_5$, 
$s_6\}$ which are never equal to zero. The one term $s_k$ which is zero 
counts as a factor equal to $1$ in the product~(\ref{eq5}).

Hence, the total number of way $N_1$ in which there could be only a 
single pair of {\em consecutive numbers} within the drawn sequence is 
then given by:

\begin{equation}
N_1=\frac{1}{84!}\cdot {5 \choose 1}\cdot \frac{d^{84}}{dx^{84}} 
\biggl[\bigg(
\frac{x^{85}-1}{x-1}\biggl)^{\bf 2}\bigg(
\frac{x^{85}-1}{x-1}-1\biggl)^{\bf 4}\biggr]
\Biggl|_{x=0}=164,007,585.
\label{eq6}
\end{equation}

Again, the factor $\frac{1}{84!}$ is the normalization factor needed to 
simplify the `spurious' numerical coefficient which originates from the 
exponent of $x^{84}$ after multiple derivative.

It should be now evident how to calculate the number of possible drawn 
sequences with {\em consecutive numbers} originating from two, three, 
four or all the five numbers $s_k$ being equal to zero. Let us call 
$N_2$, $N_3$, $N_4$, $N_5$ such numbers with respectively two, three, 
four and all the five numbers $s_k$ equal to zero. Applying the same 
reasoning as for eq.~(\ref{eq6}), one has:

\begin{equation}
N_2=\frac{1}{84!}\cdot {5 \choose 2}\cdot \frac{d^{84}}{dx^{84}}
\biggl[\bigg(\frac{x^{85}-1}{x-1}\biggl)^{\bf 2} \bigg(
\frac{x^{85}-1}{x-1}-1\biggl)^{\bf 3}\biggr]
\Biggl|_{x=0}=20,247,850.
\label{eq7}
\end{equation}

\begin{equation}
N_3=\frac{1}{84!}\cdot {5 \choose 3}\cdot \frac{d^{84}}{dx^{84}}
\biggl[\bigg(
\frac{x^{85}-1}{x-1}\biggl)^{\bf 2}\bigg(
\frac{x^{85}-1}{x-1}-1\biggl)^{\bf 2}\biggr]\Biggl|_{x=0}=987,700
\label{eq8}
\end{equation}

\begin{equation}
N_4=\frac{1}{84!}\cdot {5 \choose 4}\cdot \frac{d^{84}}{dx^{84}}
\biggl[\bigg(
\frac{x^{85}-1}{x-1}\biggl)^{\bf 2}\bigg(
\frac{x^{85}-1}{x-1}-1\biggl)\biggr]\Biggl|_{x=0}=17,850.
\label{eq9}
\end{equation}

\begin{equation}
N_5=\frac{1}{84!}\cdot {5 \choose 5}\cdot \frac{d^{84}}{dx^{84}}
\biggl[\bigg(\frac{x^{85}-1}{x-1}\biggl)^{\bf 2}\biggr]\Biggl|_{x=0}=85.
\label{eq10}
\end{equation}

The numbers $N_i$ thus count the amount of sequences with $i$ pairs of 
consecutive numbers inside (for example, in the sequence 
$\{1,3,13,14,15,87\}$, two pairs of consecutive numbers are present, 
$\{13,14\}$ and $\{14,15\}$, and hence this sequence is counted in 
$N_2$).

Summarizing, the sought number $N_c$ is thus equal to the sum:

\begin{equation}
N_c=N_1+N_2+N_3+N_4+N_5=185,261,070
\label{eq11}
\end{equation}
and the sought {\em consecutive number} probability is then:

\begin{equation}
{\cal C}(6,90)=\frac{N_c}{N_t}=\frac{185,261,070}{622,614,630}\approx 29.75\%.
\label{eq12}
\end{equation}

This outcome says that, on average, nearly one drawing every three 
results in the {\em consecutive number} phenomenon; it is definitely a 
quite common result.

A direct comparison between expected and observed {\em consecutive 
number} frequencies can be easily made taking the data in 
Table~\ref{tab1} as a statistical sample. The observed frequency amounts 
to $\frac{13}{40}\sim 32\%$, which is slightly higher than that 
predicted by eq.~(\ref{eq12}). This is probably due to the smallness of 
the statistical sample taken into account.

According to the official statistics~\cite{off} on all the drawings up 
to now (spanning from the last days of 1997 to the whole year 2009), the 
number of {\em consecutive number} occurrences $N_c$ is equal to $454$, 
while the total number of drawings $N_t$ amounts to $1507$; hence, the 
observed frequency calculated on the widest possible statistical sample 
results to be $\frac{454}{1507}\sim 30,1\%$. This figure is reasonably 
close to that predicted by eq.~(\ref{eq12}).

Another useful result that can be easily derived is the probability 
$P(n,M)$ that among a total of $M$ drawings $n$ drawings occur with {\em 
consecutive numbers}. It is not difficult to realize that $P(n,M)$ can 
be expressed as a straightforward binomial distribution with success 
probability equal to ${\cal C}(6,90)$:

\begin{equation}
P(n,M)={M\choose n}{\cal C}(6,90)^n(1-{\cal C}(6,90))^{M-n},
\label{eq13}
\end{equation}
with mean $\mu$ and standard deviation $\sigma$ equal to: 

\begin{table}
\begin{center}
{\footnotesize
\begin{tabular}{|c|c|c|}
\hline
Consecutive number categories & Observed occurrences~\cite{off} & Expected occurrences (eq.~\ref{eq13})\\
   {\footnotesize (see text for the meaning of $N_i$)}      &             &   $\mu\pm\sigma$\\
\hline  
 $N_1$  & 396 &  $\sim 397\pm 17 $\\
 $N_2$  & 53 & $\sim 49\pm 7$ \\
 $N_3$  & 5 &  $\sim 2\pm 2$\\
 $N_4$  & 0 & $\sim 0 $\\
 $N_5$  & 0 & $\sim  0$ \\
\hline
 $N_c$  & $454$ & $448\pm 18$ \\
\hline
\end{tabular}
}
\caption{Detailed comparison between observed~\cite{off} and expected 
(eq.~(\ref{eq13})) occurrences of {\em consecutive number} strings. As 
explained in Section~2, $N_i$ stands for the number of sequence with $i$ 
pairs of consecutive numbers (e.g.~a sequence like $\{1,3,13,14,15,87\}$ 
has two pairs of consecutive numbers, $\{13,14\}$ and $\{14,15\}$). The 
observed occurrences are counted among all the 1507 drawings from 
December 7, 1997, up to the whole year 2009. The expected occurrences 
are calculated with the mean and the standard deviation of the binomial 
distribution (eq.~(\ref{eq13})), with success probability equal to 
$\frac{N_i}{N_t}$.}
\end{center}
\label{tab2}
\end{table}

\begin{center}
\begin{tabular}{ccl}
$\mu$ & $=$ & $M\cdot{\cal C}(6,90)$ \\
 & & \\
$\sigma$ & $=$ & $\sqrt{M\cdot {\cal C}(6,90)\cdot(1-{\cal C}(6,90))}$ 
\end{tabular}
\end{center}

In the case of $M=1507$ (namely the number of all the {\em 
SuperEnalotto} drawings to date) the mean $\mu$ of the 
distribution~(\ref{eq13}) is equal to 448 and the standard deviation 
$\sigma$ is nearly equal to $18$. Thus, the observed number of drawings 
with {\em consecutive number} strings, $454$, is well within the 
expected one $448\pm 18$. Table~\ref{tab2} contains also a detailed 
comparison between observed~\cite{off} and expected (eq.~(\ref{eq13})) 
occurrences of {\em consecutive number} strings for each category 
$N_1$, $N_2$, $N_3$, $N_4$, $N_5$.

\section{Concluding Remarks}

In this paper a way has been suggested for calculating the probability 
of {\em consecutive number} strings within a sequence of $n$ numbers 
randomly drawn (without replacement) among the set of the first $N$ 
consecutive numbers, with $N\gg n$.

An explicit derivation has been carried out for the special case of {\em 
SuperEnalotto}, nowadays the most famous lottery in Italy, with $N=90$ 
and $n=6$. It turned out that, on average, one every three drawings 
presents one or more {\em consecutive number} strings inside: {\em a 
posteriori}, this is obviously not surprising; {\em a priori}, on the 
other hand, it admittedly appears so, at least to the author.

One reasonably expects that if the ratio $n/N$ increases, then the {\em 
consecutive number} probability will become higher. In the limit, it is 
not difficult to prove that if $n\geq \frac{N+2}{2}$, for $N$ even, or 
if $n\geq \frac{N+3}{2}$, for $N$ odd, the probability is equal to $1$. 
This comes from a direct application of the Pigeonhole Principle.

One also may think to use the results obtained in Section~2, along with 
some statistical goodness-to-fit tests to quantitatively compare 
observed and expected occurrences, as a quality control on the lottery 
operators activity. Similar kinds of check are now common against 
accounting data and tax frauds, e.g.~using the first-digit distribution 
or Benford's law~\cite{nig}.

Artifactual, man-made number sequences usually show a lower occurrence 
of {\em consecutive number} strings: commonly, people tend to severely 
underestimate the presence of {\em consecutive number} strings within 
randomly generated number sequences.

\section*{Acknowledgments}

The author acknowledges the partial support of the Italian Space Agency 
under ASI Contract N$^\circ$ 1/015/07/0. Dr.ssa Assunta Tataranni is also
kindly acknowledged for her comments and suggestions.

\end{document}